\documentclass[12pt,reqno]{amsart}

\usepackage{amsmath,amsthm,latexsym,amssymb}
\usepackage{amsmath}
\usepackage{amsfonts}
\usepackage{amssymb}
\usepackage{graphicx}
\usepackage{color}
\usepackage{cite}

\providecommand{\U}[1]{\protect\rule{.1in}{.1in}}
\numberwithin{equation}{section}

\newtheorem{theorem}{Theorem}
\newtheorem{lemma}{Lemma}

\newtheorem{remark}{Remark}

\numberwithin{theorem}{section}
\numberwithin{corollary}{section}
\numberwithin{lemma}{section}
\numberwithin{definition}{section}
\numberwithin{proposition}{section}
\numberwithin{remark}{section}

\newcommand{\medint}{-\kern  -,375cm\int}

\begin{document}
\title[ Optimal Szeg\"{o}-Weinberger type inequalities]
{ Optimal Szeg\"{o}-Weinberger type inequalities}
\author{ F. Brock$^{1}$ - F. Chiacchio$^{2}$ - G. di Blasio$^{3}$}
\thanks{}
\date{}

\begin{abstract}
Denote with $\mu _{1}(\Omega ;e^{h\left( |x|\right) })$ the first nontrivial
eigenvalue of the Neumann problem 
\begin{equation*}
\left\{ 
\begin{array}{lll}
-\text{div}\left( e^{h\left( |x|\right) }\nabla u\right) =\mu e^{h\left(
|x|\right) }u & \text{in} & \Omega \\ 
&  &  \\ 
\frac{\partial u}{\partial \nu }=0 & \text{on} & \partial \Omega ,
\end{array}
\right.
\end{equation*}
where $\Omega $ is a bounded and Lipschitz domain in $\mathbb{R}^{N}$. Under
suitable assumption on $h$ we prove that the ball centered at the origin is
the unique set maximizing $\mu _{1}(\Omega ;e^{h\left( |x|\right)})$ among
all Lipschitz bounded domains $\Omega $ of $\mathbb{R}^{N}$ of prescribed $
e^{h\left( |x|\right) }dx$-measure and symmetric about the origin. Moreover,
an example in the model case $h\left( |x|\right) =|x|^{2},$ shows that, in
general, the assumption on the symmetry of the domain cannot be dropped. In
the one-dimensional case, i.e. when $\Omega $ reduces to an interval $(a,b), 
$ we consider a wide class of weights (including both Gaussian and
anti-Gaussian). We then describe the behavior of the eigenvalue as the
interval $(a,b)$ slides along the $x$-axis keeping fixed its weighted length.

\bigskip


\end{abstract}

\maketitle

\setcounter{footnote}{1} \footnotetext{
Leipzig University, Department of Mathematics, Augustusplatz, 04109 Leipzig,
Germany, e-mail: brock@math.uni-leipzig.de}

\setcounter{footnote}{2} \footnotetext{
Dipartimento di Matematica e Applicazioni \textquotedblleft R.
Caccioppoli\textquotedblright , Universit\`{a} degli Studi di Napoli
\textquotedblleft Fe\-derico II\textquotedblright , Complesso Monte S.
Angelo, via Cintia, 80126 Napoli, Italy, e-mail: francesco.chiacchio@unina.it
}

\setcounter{footnote}{3} \footnotetext{
Dipartimento di Matematica e Fisica, Seconda Universit\`{a} degli Studi di
Napoli, via Vivaldi, 81100 Caserta, Italy, e-mail: giuseppina.diblasio@unina2.it.}

\section{Introduction}

In \cite{KS} Kornhauser and Stakgold made a famous conjecture: among all
planar simply connected domains, with fixed Lebesgue measure the first
nontrivial eigenvalue of the Neumann Laplacian achieves its maximum value if
and only if $\Omega $ is a disk. This conjecture was proved by Szeg\"{o} in 
\cite{S}. In \cite{W} Weinberger generalized this result to any bounded
smooth domain of $\mathbb{R}^{N}$. \ Adapting Weinberger arguments, similar
inequalities for spaces of constant sectional curvature have been derived
(see, e.g., \cite{AB} and \cite{cha1}). Further, it is proved in \cite{LS}
that the first nonzero Neumann eigenvalue is maximal for the equilateral
triangle among all triangles of given perimeter, and hence among all
triangles of given area.

For broad surveys of isoperimetric eigenvalue inequalities, one can consult,
for instance, the monographs of Bandle \cite{Ba}, Henrot \cite{H1} and \cite
{H2}, Kesavan \cite{Ke}, and the survey paper by Ashbaugh \cite{A}.

In this paper we derive some sharp Szeg\"{o}-Weinberger type inequalities
for the first nontrivial eigenvalue, $\mu _{1}(\Omega ;e^{h\left( |x|\right)
}),$ of the following class of problems

\begin{equation}
\left\{ 
\begin{array}{lll}
-\text{div}\left( e^{h\left( |x|\right) }\nabla u\right) =\mu e^{h\left(
|x|\right) }u & \text{in} & \Omega \\ 
&  &  \\ 
\frac{\partial u}{\partial \nu }=0 & \text{on} & \partial \Omega .
\end{array}
\right.  \label{P_1}
\end{equation}
Here and in the sequel, $\Omega $ will denote a bounded domain in $\mathbb{R}
^{N}$ with Lipschitz boundary and $\nu $ the outward normal to $\partial
\Omega $. Since the degeneracy of the operator is given in terms of the
radial function $e^{h\left( |x|\right) }$, it appears natural to let $\Omega 
$ vary in the class of sets having prescribed $\gamma _{h}$-measure, where 
\begin{equation}
d\gamma _{h}=e^{h(|x|)}dx,\text{ with }x\in \mathbb{R}^{N}.  \label{gamma_h}
\end{equation}
Recently, in \cite{CdB}, it has been proved that among all Lipschitz domains 
$\Omega $ in $\mathbb{R}^{N}$, which are symmetric about the origin and have
fixed Gaussian measure, $\mu _{1}(\Omega ;e^{-|x|^{2}/2})$ achieves its
maximum value if and only if $\Omega $ is the Euclidean ball. In the same
paper it has been shown that $\mu _{1}((a,b);e^{-t^{2}/2})$ is minimal when
the interval $(a,b)$ reduces to a half-line and maximal when it is centered
at the origin, and it is strictly monotone as $(a,b)$ slides between these
extreme positions.

The first part of the present paper deals with a class of weighted
eigenvalue problems in the form 
\begin{equation}
\left\{ 
\begin{array}{ll}
-(u^{\prime }q)^{\prime }=\mu uq & \text{ in }(a,b) \\ 
&  \\ 
u^{\prime }(a)=u^{\prime }(b)=0, & 
\end{array}
\right.  \label{P_q}
\end{equation}
where $a,b\in \mathbb{R}$ with $a<b$ and 
\begin{equation}
q(x)\in C^{2}\left( \mathbb{R}\right) ,\text{ \ }q\left( x\right) >0\text{ }
\forall x\in \mathbb{R}.  \label{q}
\end{equation}
We analyze the behavior of $\mu _{1}(\left( a,b\right) ;q)$ as the interval $
\left( a,b\right) $ slides along the $x$-axis keeping fixed its $q-$length $
\left( =\int_{a}^{b}q(x)dx\right) $. We prove that if $q(x)=q(\left\vert
x\right\vert )$ and it is decreasing on $\mathbb{R}_{+}$ then $\mu
_{1}(\left( a,b\right) ;q)$ behaves like $\mu _{1}((a,b);e^{-t^{2}/2}),$
while if $q(x)$ is increasing on $\mathbb{R}_{+}$ then $\mu _{1}(\left(
a,b\right) ;q)$ behaves in the opposite way. Note that the method used here
is different and more general than the one of \cite{CdB}. For the precise
statement see Theorem \ref{Th1 copy(1)} below. It treats the Gaussian and
anti-Gaussian-type weights in a unified way. We emphasize that, in contrast
to the case $N\geq 2$, no concavity assumptions are imposed on the weight
function. This explains the different notation in (\ref{P_1}) and
 (\ref{P_q}).

Let $c:=\int_{\mathbb{R}}q\left( t\right) dt\in \left( 0,+\infty \right] $,
and fix $d\in \left( 0,c\right) $. We define $\overset{\_}{a} >0$ such that 
\begin{equation*}
\int_{-\overset{\_}{a}}^{\overset{\_}{a}}q\left( t\right) dt=d
\end{equation*}
If $c<+\infty ,$ then we define $a_{+}$ by 
\begin{equation*}
\int_{a_{+}}^{+\infty }q\left( t\right) dt=d,
\end{equation*}
otherwise, put $a_{+}=+\infty .$

We consider the function $b=b(a)$ defined on $\left( -\infty ,a_{+}\right) $
such that 
\begin{equation*}
\int_{a}^{b\left( a\right) }q\left( t\right) =d.
\end{equation*}
Denote by $\mu _{1}\left( a\right) $ the eigenvalue $\mu _{1}\left( \left(
a,b\left( a\right) \right) ;q\right) $. Then we have the following

\vspace{0.3 cm}

\begin{theorem}
\label{Th1 copy(1)} Assume that $q$ satisfies (\ref{q}) and is even. Then

\begin{itemize}
\item[\textsl{(i)}] $\mu _{1}\left( a\right) =\mu _{1}\left( -b\left(
a\right) \right) $ for $a<a_{+},$\vspace{0.3cm}

\item[\textsl{(ii)}] if $q^{\prime }\left( x\right) \underset{(\leq )}{\geq }
0\hspace{0.1cm}$ for $x\geq 0,$ then $\dfrac{d}{da}\mu _{1}\left( a\right) 
\underset{(\leq )}{\geq }0$ for $a>-\overset{\_}{a},$\vspace{0.3cm}

\item[\textsl{(iii)}] if $q$ is not constant on 
$\left( a,b\left( a\right) \right) $ and $q^{\prime }\left( x\right) \underset{(\leq )}{\geq }0
 \hspace{0.1cm}$ for 
$x\geq 0,$ then $\dfrac{d}{da}\mu _{1}\left( a\right) 
\underset{
(<)}{>}0$ for $a>-\overset{\_}{a}.$
\end{itemize}
\end{theorem}

\bigskip

Now let us turn our attention to the $N$-dimensional problem (\ref{P_1}).

We assume that the function $h(r)$ fulfills the following set of hypotheses 
\begin{equation}
h(r)\in C^{2}\left( \left[ 0,+\infty \right[ \right) ,\text{ \ \ }h^{\prime
}(r)>-\frac{N-1}{r}\text{ \ \ }\forall r>0,\text{ }\text{\ }h^{\prime \prime
}(r)\geq 0\text{ \ \ }\forall r\geq 0.  \label{h(r)_assump}
\end{equation}

Our second main result is the following.

\begin{theorem}
\label{Main_th}Assume that assumptions (\ref{h(r)_assump}) are in force.
Then the ball centered at the origin is the unique set maximizing $\mu
_{1}(\Omega ;e^{h\left( |x|\right) })$ among all Lipschitz bounded domains $
\Omega $ of $\mathbb{R}^{N}$ of prescribed $\gamma _{h}$-measure and
symmetric about the origin. Moreover, if the assumption on the symmetry of
the domain is dropped, then, in general, the thesis does not hold true.
\end{theorem}

Finally note that Theorem \ref{Main_th} in particular applies to $\mu
_{1}(\Omega ;e^{|x|^{2}})$. Indeed, as model case we may choose $h\left(
|x|\right) =\left\vert x\right\vert ^{2}.$ Hence our result gives
information about the Neumann eigenvalues of the problem 
\begin{equation*}
-\Delta u-2x\cdot \nabla u=\mu u,
\end{equation*}
which is widely studied in literature (see, e.g., \cite{BMP} and 
\cite{BCM}). 
For related results see also \cite{BeLi, BBMP, BCKT} and \cite{DPG1, DPG2}.

The paper is organized as follows. Section 2 contains some results from the
theory of weighted rearrangements along with the definition of the suitable
Sobolev spaces naturally associated to problem (\ref{P_1}). Section 3 is
devoted to Theorem \ref{Th1 copy(1)}. By a suitable change of variable, we
firstly show that $\mu _{1}\left( \left( a,b\right) ;q\right) $ coincides
with$\ \lambda _{1}\left( \left( a,b\right) ;q^{-1}\right) $, the first
Dirichlet eigenvalue of problem (\ref{P_q}) with respect the weight $q^{-1}$
. In turn, we observe that the following problems 
\begin{equation*}
P_{q}:\left\{ 
\begin{array}{ll}
-\frac{d}{dx}\left( \frac{1}{q(x)}\frac{dv}{dx}\right) =\lambda \frac{1}{
q(x) }v(x) & \text{in }(a,b) \\ 
&  \\ 
v(a)=v(b)=0 & 
\end{array}
\right. \text{ \hspace{0.3cm} and \hspace{0.3cm} }P_{m}:\left\{ 
\begin{array}{ll}
-\frac{d^{2}w}{dy^{2}}=\lambda m(y)w(y) & \text{in }(\alpha ,\beta ) \\ 
&  \\ 
w(\alpha )=w(\beta )=0, & 
\end{array}
\right.
\end{equation*}
are isospectral, where $y=:F(x):=\int_{0}^{x}q\left( t\right) dt,\ \alpha
:=F\left( a\right) ,$ $\beta :=F\left( b\right) $ and $m$ is defined in (\ref
{m}). The advantages of studying $P_{m}$ in place of $P_{q}$ are twofold.
First, since the function $F(x)$ pushes the
measure $q(x)dx$ forward $dy$,
as long as the interval $(a,b)$ moves along the $x-$axis with fixed 
$q-$length, then $(\alpha ,\beta )$ slides along the $y-$axis keeping fixed
just its Lebesgue measure. Second, the new equation contains a weight only
in the zero order term, so it nicely behaves under reflection with respect
to $(\alpha +\beta )/2.$ These two circumstances allow to evaluate the sign
of the shape derivative of the first eigenvalue of $P_{m}$ and hence of $\mu
_{1}\left( \left( a,b\right) ;q\right) $. In Section 4 we prove Theorem \ref
{Main_th}. To this aim, we first study problem (\ref{P_1}) in the radial
case, i.e., when $\Omega =B_{R}$, the ball centered at the origin of radius $
R.$ We deduce that $\mu _{1}(B_{R};e^{h\left( |x|\right) })$ is an
eigenvalue of multiplicity $N$, and a corresponding set of linearly
independent eigenfunctions is $\left\{ w\left( \left\vert x\right\vert
\right) \dfrac{x_{i}}{\left\vert x\right\vert },\,i=1,...,N\right\} ,$ with
an appropriate function $w$. We then
define $G(r)=w(r)$ for $0\leq r\leq R$ and $G(r)=w(R)$ for $r>R$. Since the
functions $G\left( \left\vert x\right\vert \right) \dfrac{x_{i}}{\left\vert
x\right\vert }$, ($i=1,...,N$), have mean value zero, we may use them in the
variational characterization of $\mu _{1}(\Omega ;e^{h\left( |x|\right) })$.
Then the result is achieved by symmetrization arguments.

\section{Notation and preliminary results}

Now we recall a few definitions and properties about weighted rearrangement.
For exhaustive treatment on this subject we refer, e.g., to , \cite{CR}, 
\cite{Ka} and \cite{Rs}.

Throughout this paper, $B_{R}$ will denote the ball of $\mathbb{R}^{N}$
centered at the origin of radius $R.$

Let $u:x\in \Omega \rightarrow \mathbb{R}$ be a measurable function. We
denote by $m(t)$ the distribution function of $\left\vert u(x)\right\vert $
with respect to $\gamma _{h}-$measure, defined in (\ref{gamma_h}), i.e. 
\begin{equation*}
m(t)=\gamma _{h}\left( \left\{ x\in \Omega :\left\vert u(x)\right\vert
>t\right\} \right) ,\quad t\geq 0,
\end{equation*}
while the decreasing rearrangement and the increasing rearrangement of $u$
are defined respectively by 
\begin{equation*}
u^{\ast }\left( s\right) =\inf \left\{ t\geq 0:m\left( t\right) \leq
s\right\} \text{,}\quad s\in \left] 0,\gamma _{h}(\Omega )\right]
\end{equation*}
and 
\begin{equation*}
u_{\ast }\left( s\right) =u^{\ast }\left( \gamma _{h}\left( \Omega \right)
-s\right) \text{,}\quad s\in \left[ 0,\gamma _{h}(\Omega )\right[ .
\end{equation*}
Finally $u^{\bigstar }$, the $\gamma _{h}-$rearrangement of $u$, is given by 
\begin{equation*}
u^{\bigstar }\left( x\right) =u^{\bigstar }\left( |x|\right) =u^{\ast
}\left( \gamma _{h}\left( B_{\left\vert x\right\vert }\right) \right) ,\quad
x\in \Omega ^{\bigstar },
\end{equation*}
where $\Omega ^{\bigstar }$ is the ball $B_{r^{\bigstar }}$ such that $
\gamma _{h}\left( \Omega ^{\bigstar }\right) =\gamma _{h}\left(
B_{r^{\bigstar }}\right) .$ By its very definition $u^{\bigstar }$ is a
radial and radially decreasing function. Since $u$ and $u^{\bigstar }$ are $
\gamma _{h}-$equimeasurable, Cavalieri's principle ensures 
\begin{equation*}
\left\Vert u\right\Vert _{L^{p}\left( \Omega ;\gamma _{h}\right)
}=\left\Vert u^{\bigstar }\right\Vert _{L^{p}\left( \Omega ^{\bigstar
};\gamma _{h}\right) },\quad \forall p\geq 1.
\end{equation*}
We will also make use of the Hardy-Littlewood inequality, which states that 
\begin{equation}
\int_{0}^{\gamma _{h}\left( \Omega \right) }u^{\ast }\left( s\right) v_{\ast
}\left( s\right) ds\leq \int_{\Omega }\left\vert u\left( x\right) v\left(
x\right) \right\vert d\gamma _{h}\leq \int_{0}^{\gamma _{h}\left( \Omega
\right) }u^{\ast }\left( s\right) v^{\ast }\left( s\right) ds.  \label{HL}
\end{equation}

Since $\Omega $ is a bounded domain, assumptions (\ref{h(r)_assump}) ensures
that 
\begin{equation}
0<c_{1}<e^{h(\left\vert x\right\vert )}<c_{2}<+\infty ,\quad \text{in }
\Omega ,  \label{c_1_c_2}
\end{equation}
for some constants $c_{1}$ and $c_{2}.$

The natural functional space associated to problem (\ref{P_1}) is the
weighted Sobolev space defined as follows 
\begin{equation*}
H^{1}(\Omega ;\gamma _{h})=\left\{ u\in W_{\text{loc}}^{1,1}(\Omega ):\left(
u,\left\vert \nabla u\right\vert \right) \in L^{2}\left( \Omega ;\gamma
_{h}\right) \times L^{2}\left( \Omega ;\gamma _{h}\right) \right\} ,
\end{equation*}
endowed with the norm 
\begin{equation}
\left\Vert u\right\Vert _{H^{1}(\Omega ;\gamma _{h})}^{2}=\left\Vert
u\right\Vert _{L^{2}(\Omega ;\gamma _{h})}^{2}+\left\Vert \left\vert \nabla
u\right\vert \right\Vert _{L^{2}(\Omega ;\gamma _{h})}^{2}=\int_{\Omega
}u^{2}d\gamma _{h}+\int_{\Omega }\left\vert \nabla u\right\vert ^{2}d\gamma
_{h}.  \label{norm_h}
\end{equation}
By (\ref{c_1_c_2}) one immediately infers that 
\begin{equation*}
u\in H^{1}(\Omega ;\gamma _{h})\Longleftrightarrow u\in H^{1}(\Omega ),
\end{equation*}
since the norm defined in (\ref{norm_h}) is equivalent to the usual norm in $
H^{1}(\Omega )\,.$ Hence we have that $H^{1}(\Omega ;\gamma _{h})$ is
compactly embedded in $L^{2}(\Omega ;\gamma _{h})$. By standard theory on
self-adjoint compact operator, $\mu _{1}(\Omega ;e^{h\left( |x|\right) })$
admits the following well known variational characterization 
\begin{equation}
\mu _{1}(\Omega ;e^{h\left( |x|\right) })=\min \left\{ \frac{\int_{\Omega
}\left\vert \nabla v\right\vert ^{2}d\gamma _{h}}{\int_{\Omega }v^{2}d\gamma
_{h}}:v\in H^{1}(\Omega )\backslash \left\{ 0\right\} ,\text{ }\int_{\Omega
}vd\gamma _{h}=0\right\} .  \label{def_auto}
\end{equation}

\section{The one-dimensional case}

Throughout this Section we will assume that condition (\ref{q}) is fulfilled
and that 
\begin{equation*}
-\infty<a<b<+\infty .
\end{equation*}

We consider the weighted Neumann eigenvalue problem 
\begin{equation}
\left\{ 
\begin{array}{ll}
-(u^{\prime }q)^{\prime }=\mu uq & \text{ in }(a,b) \\ 
&  \\ 
u^{\prime }(a)=u^{\prime }(b)=0. & 
\end{array}
\right.   \label{1d_N}
\end{equation}
The first nontrivial eigenvalue of (\ref{1d_N}),  
$\mu_{1}\left( \left( a,b\right) ;q\right) $, clearly fulfills 
\begin{equation}
\mu _{1}\left( \left( a,b\right) ;q\right) =\min \left\{ \frac{
\int_{a}^{b}(u^{\prime })^{2}q\,dx}{\int_{a}^{b}u^{2}q\,dx}:u\in H^{1}\left(
a,b\right) \backslash \left\{ 0\right\} ,\,\int_{a}^{b}uq\,dx=0\right\} .
\label{mu_1(a,b)}
\end{equation}
Here we are interested in studying the behavior of $\mu _{1}\left( \left(
a,b\right) ;q\right) $ when the interval $(a,b)$ slides along the $x$-axis,
keeping fixed its weighted length $\int_{a}^{b}q\left( x\right) dx$. 

Further, we consider the weighted Dirichlet eigenvalue problem 
\begin{equation}
\left\{ 
\begin{array}{ll}
-\left( v^{\prime }\dfrac{1}{q}\right) ^{\prime }=\lambda v\dfrac{1}{q} & 
\text{ in }(a,b) \\ 
&  \\ 
v(a)=v(b)=0, & 
\end{array}
\right.   
\label{1d_D}
\end{equation}
and we denote its first eigenvalue by 
$\lambda _{1}\left( \left( a,b\right) ;q^{-1}\right) $. 
 In the next Lemma, we will simply write $\mu _{1}$ and $
\lambda _{1}$ for $\mu _{1}\left( \left( a,b\right) ;q\right) $ and $\lambda
_{1}\left( \left( a,b\right) ;q^{-1}\right) $, respectively.

\bigskip 

\begin{lemma}
\label{Lemma1} There holds 
\begin{equation}
    \label{m1=l1}
   \mu _1 = \lambda _1.
\end{equation}
Further, if $u_{1}$ is an eigenfunction to problem (\ref{1d_N})
corresponding to $\mu _{1}$, then the  function 
\begin{equation}
v:=u_{1}^{\prime }q  \label{v}
\end{equation}
is an eigenfunction to problem (\ref{1d_D}) with eigenvalue $\mu _{1} $. 
Finally, if $v_1 $ is an  eigenfunction to problem 
(\ref{1d_D}) corresponding to $\lambda _1 $, 
and if $x_{0}\in \left( a,b\right) $
is such that $v_{1}^{\prime }\left( x_{0}\right) =0,$ then 
\begin{equation}
u\left( x\right) :=\int_{x_{0}}^{x}v_1 \left( t\right) \frac{1}{q\left(
t\right) }dt  \label{u(x)}
\end{equation}
is an eigenfunction to problem (\ref{1d_N}) corresponding to $\lambda _{1} $.
Finally, there holds 
\begin{equation}
u\left( x\right) =-\frac{v_1 ^{\prime }\left( x\right) }{\lambda
_{1}q\left( x\right) }.  \label{u(x)secondo}
\end{equation}
\end{lemma}

\textbf{Proof.} The equation in (\ref{1d_N}) can be clearly written in the
following form 
\begin{equation*}
-u_{1}^{\prime \prime }-\frac{q^{\prime }}{q}u_{1}^{\prime }=\mu _{1}u_{1}
\end{equation*}
By differentiating we get 
\begin{equation}
-u_{1}^{\prime \prime \prime }-\frac{q^{\prime \prime }}{q}u_{1}^{\prime }+ 
\frac{q^{\prime 2}}{q^{2}}u_{1}^{\prime }-\frac{q^{\prime }}{q}u_{1}^{\prime
\prime }=\mu _{1}u_{1}^{\prime }.  \label{u'''}
\end{equation}
Plugging (\ref{v}) into (\ref{u'''}) one obtains that

\begin{equation*}
-\left( \frac{v^{\prime }}{q}\right) ^{\prime }=-\frac{v^{\prime
\prime }}{q}+\frac{v^{\prime }}{q^{2}}q^{\prime }=\mu _{1}\frac{v}{q},
\end{equation*}
with $v\left( a\right) =v\left( b\right) =0$. 
This means that $v$ is an eigenfunction to problem (\ref{1d_D})
corresponding to $\mu _{1}$, and thus
\begin{equation}
\label{lambda_1<mu_1}
\lambda_1 \leq \mu _1 .
\end{equation}
Next, let $u$ be defined by (\ref{u(x)}). Then obviously
\begin{equation}
u^{\prime }\left( a\right) =u^{\prime }\left( b\right) =0.
\label{u'(a)=u'(b)=0}
\end{equation}
By definition $v_1$ fulfills the equation in problem (\ref{1d_D}) with 
$ \lambda =\lambda_{1}$, i.e. 
\begin{equation}
\label{V1}
-\left( v_1^{\prime }\dfrac{1}{q}\right) ^{\prime }=\lambda_{1}   v_1\dfrac{1}{q}
\end{equation}
Integrating equation (\ref{V1})  from $x_{0}$ to $x$ and taking into account that $v_1 ^{\prime
}(x_{0})=0$ and (\ref{u(x)}), we obtain 
\begin{equation*}
-v_{1}^{\prime }\frac{1}{q}=\lambda _{1} u,
\end{equation*}
that is (\ref{u(x)secondo}).

Now, since $v_{1}=q u_{1}^{\prime }$, we have $v_{1}^{\prime }=q^{\prime
}u_{1}^{\prime }+qu_{1}^{\prime \prime }$. This, together with (\ref
{u(x)secondo}) and (\ref{u'(a)=u'(b)=0}), implies 
\begin{equation*}
\left\{ 
\begin{array}{ll}
-(u^{\prime }q)^{\prime }=\lambda _{1}uq & \text{ in }(a,b) \\ 
&  \\ 
u^{\prime }(a)=u^{\prime }(b)=0. & 
\end{array}
\right. 
\end{equation*}

\noindent 
Hence $u$ is an eigenfunction to problem (\ref{1d_N}) with eigenvalue $\lambda _1 $, and  
\begin{equation}
\label{lambda_1>mu_1}
\lambda_1 \geq \mu _1 .
\end{equation}
Now (\ref{lambda_1<mu_1}) and (\ref{lambda_1>mu_1}) imply (\ref{m1=l1}). The Lemma is proved. 
$\hfill \square $ \newline

\bigskip 

Now we introduce a new independent variable for problem (\ref{1d_D}), 
\begin{equation}
y=:F(x):=\int_{0}^{x}q\left( t\right) dt\text{ \ \ with \ }\alpha :=F\left(
a\right) \text{ \ and\ \ }\beta :=F\left( b\right) .  \label{y}
\end{equation}
Let us consider the following eigenvalue problem 
\begin{equation}
\left\{ 
\begin{array}{ll}
-w^{\prime \prime }=kmw & \text{ in }(\alpha ,\beta ) \\ 
&  \\ 
w(\alpha )=w(\beta )=0, & 
\end{array}
\right.   \label{1D_w}
\end{equation}
where

\begin{equation}
m\left( y\right) :=\frac{1}{\left[ q\left( F^{-1}\left( y\right) \right) 
\right] ^{2}}.  \label{m}
\end{equation}
We will denote by 
$
k_{n}((\alpha ,\beta );m)
$
the sequence of eigenvalues
to problem (\ref{1D_w}), arranged in increasing order. It is straightforward
to verify that $v_{n}$ is an eigenfunction to problem (\ref{1d_D})
corresponding to its $n$-th eigenvalue $\lambda _{n}((a,b);q^{-1}))$ if and
only if $w_{n}(y):=v_{n}(F^{-1}\left( y\right) ),$ is an eigenfunction to
the new problem (\ref{1D_w}) corresponding to its $n$-th eigenvalue $
k_{n}((\alpha ,\beta );m)$.
Indeed, since $F^{-1}(y)$ is a strictly increasing function admitting a zero 
($F^{-1}(y)=0$ if and only if $y=0$), we have that 
$w_{n}(y)$ and $v_{n}(x )$ have the same number of nodal
domains, therefore (see, e.g., \cite{C-H} Vol. 1 p. 454) there holds
\begin{equation}
\label{kn_ln}
k_{n}((\alpha ,\beta );m)=\lambda _{n}((a,b);q^{-1})), 
\hspace{0.2cm} \forall n \in \mathbb{N}.
\end{equation}

\bigskip

\begin{lemma} \label{Lemma2}
Let $m$ be defined by (\ref{m}) and even, that is $m\left(
y\right) =m\left( -y\right) $ $\forall y\in \left( -c/2,c/2\right) ,$ 
where $ c=\int_{\mathbb{R}}q(t)\,dt.$ 
Assume $\alpha +\beta >0.$ Let $w_{1}$ be the
eigenfunction to problem (\ref{1D_w}) corresponding to $k_{1}((\alpha ,
\beta );m)$, normalized so that $w_{1}(y)>0$ and 
\begin{equation*}
\int_{\alpha }^{\beta }w_{1}^{2}\left( y\right) m\left( y\right) dy=1.
\end{equation*}

\begin{itemize}
\item[(i)] If $m^{\prime }\left( y\right) \geq 0$ for $y>0,$ then 
\begin{equation}
-w_{1}^{\prime }\left( \beta \right) \geq w_{1}^{\prime }\left( \alpha
\right) .  \label{w'>}
\end{equation}

\item[(ii)] If $m^{\prime }\left( y\right) \leq 0$ for $y>0,$ then 
\begin{equation}
-w_{1}^{\prime }\left( \beta \right) \leq w_{1}^{\prime }\left( \alpha
\right) .  \label{w'<}
\end{equation}
\end{itemize}

Moreover inequalities (\ref{w'>}) and (\ref{w'<})\ are strict unless $m$ is
constant on $(\alpha ,\beta ).$
\end{lemma}

\bigskip

\textbf{Proof. }Setting $\overline{w}\left( y\right) :=w_{1}\left( \alpha
+\beta -y\right) ,$ it follows that 
\begin{equation*}
\left\{ 
\begin{array}{lll}
-\overline{w}^{\prime \prime }=k_{1}((\alpha ,\beta );m)\text{ }\overline{w}
(y)m\left( \alpha +\beta -y\right) & \text{on} & \left( \alpha ,\beta \right)
\\ 
&  &  \\ 
\overline{w}\left( \alpha \right) =\overline{w}\left( \beta \right) =0, &  & 
\end{array}
\right.
\end{equation*}
and 
\begin{equation*}
w_{1}\left( \frac{\alpha +\beta }{2}\right) =\overline{w}\left( \frac{\alpha
+\beta }{2}\right) .
\end{equation*}
Further, setting $W:=\overline{w}-w_{1}$, we find 
\begin{equation}
\left\{ 
\begin{array}{lll}
-W^{\prime \prime }=k_{1}((\alpha ,\beta );m)\text{ }m\left( y\right) W+g & 
\text{on} & \left( \alpha ,\frac{\alpha +\beta }{2}\right) \\ 
&  &  \\ 
W\left( \alpha \right) =W\left( \frac{\alpha +\beta }{2}\right) =0 &  & 
\end{array}
\right.  \label{W}
\end{equation}
where 
\begin{equation*}
g\left( y\right) :=k_{1}((\alpha ,\beta );m)\overline{w}\left( y\right)
\left( m\left( \alpha +\beta -y\right) -m\left( y\right) \right) .
\end{equation*}

\noindent By the well known monotonicity property of the Dirichlet eigenvalues with respect to domains inclusion, we have that
$$
k_{1}\left( \left( \alpha ,\frac{\alpha +\beta }{2}
\right) ;m\right) >k_{1}((\alpha ,\beta );m).
$$ 
From the maximum principle we
obtain that the solution $W$ of the boundary value problem (\ref{W}), with $
g $ given, is unique. Moreover, $W\equiv 0$ for $g\equiv 0,$ $W>0$ for 
$g\geq 0,g \not\equiv 0$ and $W<0$ for $g\leq 0, \hspace{0.1 cm} g\not\equiv 0.$ 
Hence, using the strong maximum principle, we get

\begin{eqnarray*}
0 &=&W^{\prime }\left( \alpha \right) =-w_{1}^{\prime }\left( \alpha \right)
-w_{1}^{\prime }\left( \beta \right) \text{ \ for \ }g\equiv 0, \\
0 &<&W^{\prime }\left( \alpha \right) =-w_{1}^{\prime }\left( \alpha \right)
-w_{1}^{\prime }\left( \beta \right) \text{ \ for  }g\geq 0,\text{ }g \not\equiv 0
\text{ \ on  }\left( \alpha ,\frac{\alpha +\beta }{2}\right) , \\
0 &>&W^{\prime }\left( \alpha \right) =-w_{1}^{\prime }\left( \alpha \right)
-w_{1}^{\prime }\left( \beta \right) \text{ \ for \ }g\leq 0,\text{ }g \not\equiv 0
\text{ \ on \ }\left( \alpha ,\frac{\alpha +\beta }{2}\right) .
\end{eqnarray*}

Note that $g\geq 0$ $(g\leq 0)$ on $\left( \alpha ,\frac{\alpha +\beta }{2}
\right) $ if $m^{\prime }\left( y\right) \geq 0$ $\left( m^{\prime }\left(
y\right) \leq 0\right) $ for $y>0,$ and \thinspace $g\equiv 0$ is possible
in either case if $m=$\textit{const. }on\textit{\ }$\left( \alpha ,\beta
\right) $. This completes the proof of Lemma.$\hfill \square $

\bigskip

Let $\alpha ,\beta \in \mathbb{R}$  be defined in (\ref{y}),  let $t>0,$ be such that  $-c/2<\alpha <\beta <\beta
+t<c/2, $ 
where $ c=\int_{\mathbb{R}}q(t)\,dt.$
We consider the following family of eigenvalue problems

\begin{equation}
\left\{ 
\begin{array}{lll}
-\widetilde{w}_{yy}(y,t)=k(t)m(y)\widetilde{w}(y,t) & \text{on} & \left(
\alpha +t,\beta +t\right) \\ 
&  &  \\ 
\widetilde{w}\left( \alpha +t,t\right) =\widetilde{w}\left( \beta
+t,t\right) =0. &  & 
\end{array}
\right.  \label{w}
\end{equation}

\noindent Let 
\begin{equation}
k_{1}(t):=k_{1}\left( \left( \alpha +t,\beta +t\right) ;m\right)
\label{k(t)}
\end{equation}
be the first eigenvalue to problem (\ref{w}),
and denote with $\widetilde{w}(y,t)$ the corresponding eigenfunction,
normalized so that $\widetilde{w}(y,t)>0$ and 
\begin{equation}
\int_{\alpha +t}^{\beta +t}\widetilde{w}\left( y,t\right) ^{2}m\left(
y\right) dy=1.  \label{normalized}
\end{equation}
Clearly 
\begin{equation*}
\widetilde{w}(y,0)=w_{1}(y),
\end{equation*}
where $w_{1}(y)$ is the function defined in Lemma \ref{Lemma2}.

\begin{lemma}
\label{lamda'}Let $k_{1}(t)$ be defined by (\ref{k(t)}). It holds that

\begin{equation}
k_{1}^{\prime }(0)=-\left( \widetilde{w}_{y}(\beta ,0)\right) ^{2}+\left( 
\widetilde{w}_{y}(\alpha ,0)\right) ^{2}=-\left( w_{1}^{\prime }(\beta
)\right) ^{2}+\left( w_{1}^{\prime }(\alpha )\right) ^{2}.  \label{formula l}
\end{equation}
\end{lemma}

\textbf{Proof.} By (\ref{normalized}) we have 
\begin{equation*}
k_{1}\left( t\right) =\int_{\alpha +t}^{\beta +t}\widetilde{w}_{y}\left(
y,t\right) ^{2}dy,
\end{equation*}
by differentiating we find 
\begin{equation}
k_{1}^{\prime }(0)=2\int_{\alpha }^{\beta }\widetilde{w}_{y}(y,0)\widetilde{
w }_{ty}(y,0)dy+\left( \widetilde{w}_{y}(\beta ,0)\right) ^{2}-\left( 
\widetilde{w}_{y}(\alpha ,0)\right) ^{2}.  \label{l=}
\end{equation}
Moreover, from (\ref{normalized}) it follows that 
\begin{equation}
\int_{\alpha }^{\beta }\widetilde{w}\left( y,t\right) \widetilde{w}
_{t}(y,0)m\left( y\right) dy=0.  \label{=0}
\end{equation}
Finally, a differentiation with respect to $t$ of the equation in problem ( 
\ref{w}) gives 
\begin{equation*}
-\widetilde{w}_{tyy}(y,0)=k_{1}^{\prime }(0)m(y)\widetilde{w}
(y,0)+k_{1}(0)m(y)\widetilde{w}_{t}(y,0).
\end{equation*}
Multiplying the above equation with $\widetilde{w}(y,0)$ and integrating
yields 
\begin{equation*}
\int_{\alpha }^{\beta }\widetilde{w}_{y}(y,0)\widetilde{w}
_{ty}(y,0)dy=k_{1}^{\prime }(0)\int_{\alpha }^{\beta }m(y)\widetilde{w}
(y,0)^{2}dy+k_{1}(0)\int_{\alpha }^{\beta }m(y)\widetilde{w}_{t}(y,0) 
\widetilde{w}(y,0)dy.
\end{equation*}
Using (\ref{normalized}), (\ref{=0}) and integrating by parts, we have 
\begin{equation*}
\int_{\alpha }^{\beta }\widetilde{w}_{y}(y,0)\widetilde{w}
_{ty}(y,0)dy=k_{1}^{\prime }(0),
\end{equation*}
this, together with (\ref{l=}), gives the claim (\ref{formula l}).$\hfill
\square $

\bigskip

\begin{lemma}
\label{Lemma3} Let $\alpha +\beta >0$ and $m(y)$ be an even function, and
let $k_{1}(t)$ be given by (\ref{k(t)}).

\begin{itemize}
\item[(i)] If $m^{\prime }\left( y\right) \geq 0$ for $y>0,$ then 
\begin{equation}
k_{1}^{\prime }(0)\leq 0.  \label{l<0}
\end{equation}

\item[(ii)] If $m^{\prime }\left( y\right) \leq 0$ for $y>0,$ then 
\begin{equation}
k_{1}^{\prime }(0)\geq 0.  \label{l>0}
\end{equation}
\end{itemize}

Moreover, the inequalities (\ref{l<0}), (\ref{l>0}) are strict unless $m$ is
constant on $\left( \alpha ,\beta \right) .$
\end{lemma}

\textbf{Proof.} Lemma \ref{lamda'} tells us that 
\begin{equation}
k_{1}^{\prime }(0)=\left( w_{1}^{\prime }(\alpha )+w_{1}^{\prime }(\beta
)\right) \left( w_{1}^{\prime }(\alpha )-w_{1}^{\prime }(\beta )\right) .
\label{L1}
\end{equation}
Since we are assuming that $w_{1}(y)>0$ in $(\alpha ,\beta ),$ it holds that 
\begin{equation}
w_{1}^{\prime }(\alpha )>0>w_{1}^{\prime }(\beta ).  \label{L2}
\end{equation}
The thesis follows from (\ref{L1}), (\ref{L2}) and Lemma \ref{Lemma2}.$
\hfill \square $

\vspace{.3cm}

Now let $\overset{\_}{a},$ $a_{+}$, $c$ and $d$ the constants defined in the
Introduction. Recalling that the function $b=b(a)$ is defined on $\left(
-\infty ,a_{+}\right) $ by the identity 
\begin{equation}
\int_{a}^{b\left( a\right) }q\left( t\right) =d.  \label{Final_d}
\end{equation}
In the $y$ variable (see (\ref{y})),  the above condition obviously becomes  
\begin{equation*}
\beta =\beta (\alpha )=\alpha +d.
\end{equation*}

\noindent Now we are in position to prove the main result of this Section.
\bigskip

\textbf{\ Proof of Theorem \ref{Th1 copy(1)}. } Since $q(x)$ is an even
function, we have that $b(-b(a))=a$ and, hence, \textsl{(i) }follows.
Using (\ref{m}) we have that $q^{\prime }\left( x\right) \geq 0$ $\forall
x\geq 0$ if and only if $m^{\prime }\left( y\right) \leq 0$ for $y\in \left(
0,c/2\right) .$ Furthermore, it is straightforward to check that 
\begin{equation*}
\alpha +\beta (\alpha )>0\Leftrightarrow a+b(a):=F^{-1}(\alpha
)+F^{-1}(\alpha +d)>0\Leftrightarrow a>-\overset{\_}{a}.
\end{equation*}
From (\ref{kn_ln}), for $n=1$, and  (\ref{m1=l1}) we deduce  that
\begin{equation*}
k_{1} \left(  ( \alpha , \alpha +d  ); m \right)
= \lambda_{1}  \left( (a,b(a) );q^{-1} \right)
=\mu _{1}\left( (a,b(a) );q \right)
=: \mu _{1} ( a ).
\end{equation*}
The above chain of equalities immediately implies that 
\begin{equation*}
\dfrac   {d}{d \alpha}k_{1} \left(  ( \alpha , \alpha +d  ); m \right) =\frac{1}{q(a)}\dfrac{d}{da}\mu _{1}(a),
\end{equation*}
which in turn, thanks to Lemma \ref{Lemma3},  yields \textsl{(ii)} and \textsl{(iii)}.
  $ \hfill \square $

\bigskip

Lemma \ref{Lemma2} also allows to obtain qualitative properties about the
first nontrivial eigenfunction to problem (\ref{1d_N}).

\bigskip

\begin{lemma}
Let $q$ be even and $u_{1}$ an eigenfunction to problem (\ref{1d_N}) with
eigenvalue $\mu _{1}$ and $a+b>0.$

\begin{itemize}
\item[(i)] If $q^{\prime }\left( x\right) \geq 0$ for $x>0,$ then 
\begin{equation}
\left\vert u_{1}\left( a\right) \right\vert \geq \left\vert u_{1}\left(
b\right) \right\vert .  \label{u_1>}
\end{equation}

\item[(ii)] If $q^{\prime }\left( x\right) \leq 0$ for $x>0,$ then 
\begin{equation}
\left\vert u_{1}\left( a\right) \right\vert \leq \left\vert u_{1}\left(
b\right) \right\vert .  \label{u_1<}
\end{equation}
\end{itemize}

Moreover, inequalities (\ref{u_1>}) and (\ref{u_1<}) are strict if $q$ is
not constant on $\left( a,b\right) .$
\end{lemma}

\bigskip

\textbf{Proof.} By Lemma \ref{Lemma1}, we may assume that 
\begin{equation}
u_{1}\left( x\right) =-\frac{v_{1}^{\prime }\left( x\right) }{\lambda
_{1}\left( \left( a,b\right) ;q^{-1}\right) q\left( x\right) },  \label{ff}
\end{equation}
where $v_{1}$ is an eigenfunction to problem (\ref{1d_D}) 
corresponding to $\lambda _{1}$.
Since $w_{1}\left( F(x)\right) =Cv_{1}\left( x\right) $ for some constant $
C\neq 0,$ where $w_{1}\left( y\right) $ is the function defined in Lemma \ref
{Lemma2}, identity (\ref{ff}) becomes 
\begin{equation*}
u_{1}\left( F^{-1}(y)\right) =-\frac{1}{C}\frac{w_{1}^{\prime }\left(
y\right) }{\lambda _{1}\left( \left( a,b\right) ;q^{-1}\right) }=-\frac{1}{C}
\frac{w_{1}^{\prime }\left( y\right) }{k_{1}\left( \left( \alpha ,\beta
\right) ;m\right) }.
\end{equation*}
Now the assertions follow from Lemma \ref{Lemma2}.$\hfill \square $

\bigskip

\section{The $N-$dimensional case}

Let us consider the problem (\ref{P_1}) in $B_{R}$, the ball centered at the
origin with radius $R$, i.e.

\begin{equation}
\left\{ 
\begin{array}{lll}
-\text{div}\left( e^{h}\nabla u\right) =\mu e^{h}u & \text{in} & B_{R} \\ 
&  &  \\ 
\dfrac{\partial u}{\partial \nu }=0 & \text{on} & \partial B_{R}.
\end{array}
\right.  \label{P_Ball}
\end{equation}
The equation in (\ref{P_Ball}) can be rewritten, using polar coordinates, as 
\begin{equation}
\frac{1}{r^{N-1}}\frac{\partial }{\partial r}\left( r^{N-1}\frac{\partial u}{
\partial r}\right) +\frac{1}{r^{2}}\Delta _{\mathbb{S}^{N-1}}\left( u| 
\mathbb{S}_{r}^{N-1}\right) +h^{\prime }(r)\frac{\partial u}{\partial r}+\mu
u=0,  \label{laplaciano polare}
\end{equation}
where $\mathbb{S}_{r}^{N-1}$ is the sphere of radius $r$ in $\mathbb{R}^{N},$
\textsl{\ }$u|\mathbb{S}_{r}^{N-1}$\textsl{\ }is the restriction of $u$ on $
\mathbb{S}_{r}^{N-1}$ and finally $\Delta _{\mathbb{S}^{N-1}}\left( u| 
\mathbb{S}_{r}^{N-1}\right) $ is the standard Laplace-Beltrami operator
relative to the manifold $\mathbb{S}_{r}^{N-1}.$

Looking for separated solutions $u\left( x\right) =Y\left( \theta \right)
f\left( r\right) $\ of equation (\ref{laplaciano polare}), where $\theta $
belongs to $\mathbb{S}_{1}^{N-1},$ we find 
\begin{equation*}
Y\frac{1}{r^{N-1}}\left( r^{N-1}f^{\prime }\right) ^{\prime }+\Delta _{ 
\mathbb{S}^{N-1}}Y\frac{f}{r^{2}}+Yh^{\prime }(r)f^{\prime }+\mu Yf=0,
\end{equation*}
and hence 
\begin{equation}
\frac{1}{r^{N-3}f}\left( r^{N-1}f^{\prime }\right) ^{\prime }+r^{2}h^{\prime
}(r)\frac{f^{\prime }}{f}+\mu r^{2}=-\frac{\Delta _{\mathbb{S}^{N-1}}Y}{Y}= 
\overset{\_}{k}.  \label{eq separate}
\end{equation}
As well known, see, e.g., \cite{Mu} and \cite{Cha}, the last equality is
fulfilled if and only if 
\begin{equation}
\overset{\_}{k}=k\left( k+N-2\right) \text{ \ \ with }k=\mathbb{N}\cup
\left\{ 0\right\} .  \label{k_bar}
\end{equation}
Multiplying the left hand side of equation (\ref{eq separate}) by $\dfrac{f}{
r^{2}},$ we get 
\begin{equation*}
f^{\prime \prime }+f^{\prime }\left( \frac{N-1}{r}+h^{\prime }(r)\right)
+\mu f-k\left( k+N-2\right) \dfrac{f}{r^{2}}=0\text{ \ in \ \ }\left(
0,R\right) .
\end{equation*}
Let us denote by $f_{k},$ $Y_{k}$ the solutions of (\ref{eq separate}) with $
k=\overset{\_}{k}$ defined in (\ref{k_bar}).

The eigenfunctions are either purely radial 
\begin{equation}
u_{i}\left( r\right) =f_{0}\left( \mu _{i};r\right) ,\text{ if }k=0,
\label{Rad_eigen}
\end{equation}
or have the form

\begin{equation}
u_{i}\left( r,\theta \right) =f_{k}\left( \mu _{i};r\right) Y_{k}\left(
\theta \right) ,\ \text{if}\ k\in \mathbb{N}.  \label{Ang_eigen}
\end{equation}
The functions $f_{k},$ with $k\in \mathbb{N}\cup \left\{ 0\right\} ,$ \
clearly satisfy 
\begin{equation}
\left\{ \left. 
\begin{array}{l}
f_{k}^{\prime \prime }+f_{k}^{\prime }\left( \dfrac{N-1}{r}+h^{\prime
}(r)\right) +\mu _{i}f_{k}-k\left( k+N-2\right) \dfrac{f_{k}}{r^{2}}=0\text{
\ \ in \ \ \ }\left( 0,R\right) \\ 
\\ 
f_{k}\left( 0\right) =0,\text{ \ }f_{k}^{\prime }\left( R\right) =0.
\end{array}
\right. \right.  \label{eq_for_f_k}
\end{equation}

In the sequel we will denote by $\tau _{n}(R),$ with $n\in \mathbb{N}\cup
\left\{ 0\right\} $, the increasing sequence of eigenvalues of 
(\ref{P_Ball}) whose corresponding eigenfunctions are purely radial, i.e. in the form
 (\ref{Rad_eigen}) or equivalently solutions to problem (\ref{eq_for_f_k})
with $k=0$. Clearly in this case the first eigenfunction is constant and the
corresponding eigenvalue $\tau _{0}(R)$ is trivially zero. We will denote by 
$\nu _{n}(R)$, with $n\in \mathbb{N}$, the remaining eigenvalues of (\ref
{P_Ball}), arranged in increasing order.

\bigskip

\begin{lemma}
\label{Lemma}If the function $h(r)$ fulfills assumptions (\ref{h(r)_assump})
then 
\begin{equation}
\nu _{1}(R)<\tau _{1}(R),\text{ \quad }\forall R>0.  \label{vu1<tau2}
\end{equation}
\end{lemma}

\textbf{Proof. }
We recall that $\tau _{1}=\tau _{1}(R)$ is the first
nontrivial eigenvalue of 
\begin{equation}
\left\{ 
\begin{array}{ll}
g^{\prime \prime }+g^{\prime }\left( \dfrac{N-1}{r}+h^{\prime }(r)\right)
+\tau g=0 & \mbox{ in }(0,R) \\ 
&  \\ 
g^{\prime }\left( 0\right) =g^{\prime }\left( R\right) =0, & 
\end{array}
\right.  \label{prob_2}
\end{equation}
and $\nu _{1}=\nu _{1}(R)$ is the first eigenvalue of 
\begin{equation}
\left\{ 
\begin{array}{ll}
w^{\prime \prime }+\left( \dfrac{N-1}{r}+h^{\prime }(r)\right) w^{\prime
}+\nu w-\dfrac{N-1}{r^{2}}w=0 & \mbox{ in }(0,R) \\ 
&  \\ 
w\left( 0\right) =w^{\prime }\left( R\right) =0. & 
\end{array}
\right.  \label{prob_1}
\end{equation}
First of all\textbf{\ }we observe that the first eigenfunction $w_{1}$ of ( 
\ref{prob_1}) does not change its sign in $\left( 0,R\right) $, thus we can
assume that $w_{1}>0$ in $\left( 0,R\right) .$

Moreover $w_{1}^{\prime }>0$ in $\left( 0,R\right) .$ Indeed, assume, by
contradiction, that we can find two values $r_{1}$, $r_{2},$ with $
r_{1}<r_{2},$ such that $w_{1}^{\prime \prime }\left( r_{1}\right) \leq 0,$ $
w_{1}^{\prime }\left( r_{1}\right) =0$ and $w_{1}^{\prime \prime }\left(
r_{2}\right) \geq 0,$ $w_{1}^{\prime }\left( r_{2}\right) =0.$ By evaluating
the equation in (\ref{prob_1}) 
\begin{equation*}
\frac{w_{1}^{\prime \prime }}{w_{1}}+\frac{w_{1}^{\prime }}{w_{1}}\left( 
\frac{N-1}{r}+h^{\prime }(r)\right) +\nu _{1}-\frac{N-1}{r^{2}}=0
\end{equation*}
at $r_{1}$ and $r_{2},$ we get 
\begin{equation*}
\nu _{1}-\frac{N-1}{r_{2}^{2}}\leq 0\text{ \ and \ \ }\nu _{1}-\frac{N-1}{
r_{1}^{2}}\geq 0,
\end{equation*}
which means $r_{1}\geq r_{2}$ and this is a contradiction.

On the other hand, the first nontrivial eigenfunction of problem (\ref
{prob_2}), $g_{1}=g_{1}(r)$, has mean value zero i.e. 
\begin{equation*}
\int_{B_{R}}g_{1}e^{h(\left\vert x\right\vert )}dx=N\omega
_{N}\int_{0}^{R}g_{1}(r)e^{h(r)}r^{N-1}dr=0,
\end{equation*}
where, here and in the sequel, $\omega _{N}$ denotes the Lebesgue measure of
the unit ball in $\mathbb{R}^{N}$.

This implies that $g_{1}(r)$ must change its sign in $\left( 0,R\right) $.
Let us suppose $g_{1}(r)>0$ in $\left( 0,r_{0}\right) $ and $g_{1}\left(
r_{0}\right) =0.$ We observe that $g_{1}^{\prime }(r)<0$ in $\left(
0,R\right) $ and in particular 
\begin{equation}
g_{1}^{\prime }(r_{0})<0.  \label{eq_for_g'}
\end{equation}
Therefore evaluating the equation of problem (\ref{prob_2}) at $r_{0},$ we
have 
\begin{equation}
g_{1}^{\prime \prime }\left( r_{0}\right) +g_{1}^{\prime }\left(
r_{0}\right) \left( \frac{N-1}{r_{0}}+h^{\prime }(r_{0})\right) =0
\label{eq_for_g''}
\end{equation}
and by the assumption on $h^{\prime },$ see (\ref{h(r)_assump}), it follows
that 
\begin{equation}
g_{1}^{\prime \prime }\left( r_{0}\right) >0.  \label{g''(r_0)>0}
\end{equation}
Moreover if we set $\psi =g_{1}^{\prime },$ then problem (\ref{prob_2})
becomes 
\begin{equation}
\left\{ 
\begin{array}{ll}
\psi ^{\prime \prime }+\psi ^{\prime }\left( \dfrac{N-1}{r}+h^{\prime
}(r)\right) +\psi \left( -\dfrac{N-1}{r^{2}}+h^{\prime \prime }(r)\right)
+\tau _{1}\psi =0 & \mbox{ in }\left( 0,R\right) \\ 
&  \\ 
\psi \left( 0\right) =\psi \left( R\right) =0. & 
\end{array}
\right.  \label{A}
\end{equation}
Further since we are assuming, see (\ref{h(r)_assump}), that $h^{\prime
\prime }\geq 0,$ from (\ref{prob_1}) we have that 
\begin{equation}
\left\{ 
\begin{array}{ll}
w^{\prime \prime }+\left( \dfrac{N-1}{r}+h^{\prime }(r)\right) w^{\prime
}+(\nu _{1}+h^{\prime \prime }(r))w-\dfrac{N-1}{r^{2}}w\geq 0 & \mbox{ in }
(0,R) \\ 
&  \\ 
w\left( 0\right) =w^{\prime }\left( R\right) =0. & 
\end{array}
\right.  \label{B}
\end{equation}
Now we multiply the equation in (\ref{A}) by $r^{N-1}e^{h(r)}w_{1}$ and the
equation in (\ref{B}) by $r^{N-1}e^{h(r)}\psi ,$ respectively, and, finally,
subtracting, leads to 
\begin{gather*}
r^{N-1}e^{h(r)}\left( w_{1}^{\prime \prime }\psi -w_{1}\psi ^{\prime \prime
}\right) +r^{N-1}e^{h(r)}\left( \frac{N-1}{r}+h^{\prime }(r)\right) \left(
w_{1}^{\prime }\psi -w_{1}\psi ^{\prime }\right) + \\
+(\nu _{1}-\tau _{1})r^{N-1}e^{h(r)}w_{1}\psi \geq 0\text{ \ \ in \ \ }
(0,r_{0}).
\end{gather*}
Integrating the above inequality on $(0,r_{0}),$ we get 
\begin{gather}
(\nu _{1}-\tau _{1})\int_{0}^{r_{0}}w_{1}\psi r^{N-1}e^{h(r)}dr\geq 
\label{V-T>=} \\
\int_{0}^{r_{0}}\left[ \left( w_{1}\psi ^{\prime \prime }-w_{1}^{\prime
\prime }\psi \right) +\left( \frac{N-1}{r}+h^{\prime }(r)\right) \left(
w_{1}\psi ^{\prime }-w_{1}^{\prime }\psi \right) \right] e^{h(r)}r^{N-1}dr. 
\notag
\end{gather}
Now we claim that 
\begin{equation*}
\int_{0}^{r_{0}}\left[ \left( w_{1}\psi ^{\prime \prime }-w_{1}^{\prime
\prime }\psi \right) +\left( \frac{N-1}{r}+h^{\prime }(r) \right) \left( w_{1}\psi
^{\prime }-w_{1}^{\prime }\psi \right) \right] e^{h(r)}r^{N-1}dr>0
\end{equation*}
To this aim we first note that 
\begin{eqnarray}
\int_{0}^{r_{0}}\psi ^{\prime \prime }w_{1}r^{N-1}e^{h\left( r\right) }dr
&=&r_{0}^{N-1}\psi ^{\prime }(r_{0})w_{1}(r_{0})e^{h\left( r_{0}\right) }
\label{eq_1_r0} \\
&&-\int_{0}^{r_{0}}\psi ^{\prime }\left( h^{\prime }(r)w_{1}+w_{1}^{\prime
}+ \frac{N-1}{r}w_{1}\right) e^{h\left( r\right) }r^{N-1}dr  \notag
\end{eqnarray}
and 
\begin{eqnarray}
\int_{0}^{r_{0}}w_{1}^{\prime \prime }\psi r^{N-1}e^{h\left( r\right) }dr
&=&r_{0}^{N-1}\psi (r_{0})w_{1}^{\prime }(r_{0})e^{h\left( r_{0}\right) }
\label{eq_2_r0} \\
&&-\int_{0}^{r_{0}}w_{1}^{\prime }\left( h^{\prime }(r)\psi +\psi ^{\prime
}+ \frac{N-1}{r}\psi \right) e^{h\left( r\right) }r^{N-1}dr.  \notag
\end{eqnarray}
Recalling that $\psi ^{\prime }(r_{0})=g_{1}^{\prime \prime }(r_{0})>0,$ see
(\ref{g''(r_0)>0}), and $\psi (r_{0})=g_{1}^{\prime }(r_{0})<0,$ see (\ref
{eq_for_g'}), we have 
\begin{equation}
r_{0}^{N-1}e^{h(r_{0})}\left( \psi ^{\prime }(r_{0})w_{1}(r_{0})-\psi
(r_{0})w_{1}^{\prime }(r_{0})\right) >0.  \label{eq_3_r0}
\end{equation}
Hence, subtracting equations (\ref{eq_1_r0}) and (\ref{eq_2_r0}), taking
into account of (\ref{eq_3_r0}), from (\ref{V-T>=}) we get 
\begin{equation*}
(\nu _{1}-\tau _{1})\int_{0}^{r_{0}}w_{1}\psi e^{h(r)}r^{N-1}dr>0.
\end{equation*}
Finally, since $\psi (r)=g_{1}^{\prime }(r)<0$ in $(0,R)$ and $w_{1}(r)>0$
in $(0,R),$ we have that 
\begin{equation*}
\int_{0}^{r_{0}}w_{1}\psi e^{h(r)}r^{N-1}dr<0,
\end{equation*}
and therefore $(\nu _{1}-\tau _{1})$ must be negative. The Lemma is so
proved.$\hfill \square $

\bigskip

From Lemma \ref{Lemma} we clearly have 
\begin{equation}
\mu _{1}(B_{R};e^{h(\left\vert x\right\vert )})=\frac{\displaystyle 
\int_{B_{R}}\left( \left( w^{\prime }\left( \left\vert x\right\vert \right)
\right) ^{2}+\dfrac{N-1}{\left\vert x\right\vert ^{2}}w\left( \left\vert
x\right\vert \right) ^{2}\right) d\gamma _{h}}{\displaystyle 
\int_{B_{R}}w(\left\vert x\right\vert )^{2}d\gamma _{h}},\text{ }\forall R>0.
\label{eig_BR}
\end{equation}
Now we are in position the prove the main result of this Section.

\bigskip

\textbf{Proof of Theorem \ref{Main_th}} 
\ Recall that $\Omega ^{\bigstar }$
is the ball $B_{r^{\bigstar }}$ such that $\gamma _{h}\left( \Omega
^{\bigstar }\right) =\gamma _{h}\left( B_{r^{\bigstar }}\right) .$ We define 
\begin{equation}
G(r)=\left\{ 
\begin{array}{ll}
w(r) & \mbox{ for }0<r<r^{\bigstar } \\ 
&  \\ 
w(r^{\bigstar }) & \mbox{ for }r\geq r^{\bigstar },
\end{array}
\right.   \label{def G}
\end{equation}
where $w$ is the solution to problem (\ref{prob_1}) satisfying (\ref{eig_BR}). By
the results stated above the function $G$ is nondecreasing and nonnegative.
We introduce the functions 
\begin{equation*}
P_{i}(x)=G(\left\vert x\right\vert )\frac{x_{i}}{\left\vert x\right\vert } 
\text{ \ \ for \ \ }1\leq i\leq N.
\end{equation*}
The assumption on the symmetry of $\Omega $ guarantees 
\begin{equation}
\int_{\Omega }P_{i}(x)d\gamma _{h}=0,\text{ \ \ }\forall i=1,...,N.
\label{Orth_P_i}
\end{equation}
Hence each function $P_{i}$ is admissible in the variational formulation of $
\mu _{1}(\Omega ;e^{h\left( |x|\right) }),$ i.e. (\ref{def_auto}).

Since 
\begin{equation*}
\frac{\partial P_{i}}{\partial x_{j}}=G^{\prime }(\left\vert x\right\vert ) 
\frac{x_{i}x_{j}}{\left\vert x\right\vert ^{2}}-G(\left\vert x\right\vert ) 
\frac{x_{i}x_{j}}{\left\vert x\right\vert ^{3}}+\delta _{ij}\frac{
G(\left\vert x\right\vert )}{\left\vert x\right\vert },
\end{equation*}
where $\delta _{ij}$ is the Kronecker symbol. Using $P_{i}$ as trial
functions for $\mu _{1}(\Omega ;e^{h\left( |x|\right) })$ we get 
\begin{equation}
\mu _{1}(\Omega ;e^{h\left( |x|\right) })\leq \frac{\displaystyle 
\sum\limits_{i=1}^{N}\int_{\Omega }\sum\limits_{j=1}^{N}\left( \frac{
\partial P_{i}}{\partial x_{j}}\right) ^{2}d\gamma _{h}}{\displaystyle 
\sum\limits_{i=1}^{N}\int_{\Omega }P_{i}^{2}d\gamma _{h}}=\frac{ 
\displaystyle \int_{\Omega }N(\left\vert x\right\vert )d\gamma _{h}}{ 
\displaystyle \int_{\Omega }D(\left\vert x\right\vert )d\gamma _{h}}.
\label{mu_1_i}
\end{equation}
where 
\begin{equation*}
N(r)=\left( G^{\prime }\left( r\right) \right) ^{2}+\dfrac{N-1}{r^{2}}
G^{2}\left( r\right)
\end{equation*}
and 
\begin{equation*}
D(r)=G^{2}\left( r\right) .
\end{equation*}
We claim that 
\begin{equation*}
\frac{d}{dr}N(r)<0.
\end{equation*}
Taking into account of the definition (\ref{def G}) of $G$, and the
differential equation in (\ref{prob_1}), with $\nu =\nu _{1},$ we have 
\begin{equation*}
\frac{d}{dr}N(r)=\left\{ 
\begin{array}{ccc}
-2\left[ \nu _{1}ww^{\prime }+\left( w^{\prime }\right) ^{2}r+\dfrac{N-1}{
r^{3}}\left( rw^{\prime }+w\right) ^{2} \right] & \text{if} & 0<r<r^{\bigstar
} \\ 
&  &  \\ 
-2\left( N-1\right) \dfrac{w^{2}(r^{\bigstar })}{\left( r^{\bigstar }\right)
^{3}} & \text{if} & r\geq r^{\bigstar }
\end{array}
\right.
\end{equation*}
Now we claim that 
\begin{equation}
\int_{\Omega }N(\left\vert x\right\vert )d\gamma _{h}\leq \int_{\Omega
^{\bigstar }}N(\left\vert x\right\vert )d\gamma _{h}.  \label{dis_BR_1}
\end{equation}
Hardy-Littlewood inequality (\ref{HL}) ensures 
\begin{equation}
\int_{\Omega }N(\left\vert x\right\vert )d\gamma _{h}\leq \int_{0}^{\gamma
_{h}(\Omega )}N^{\ast }(s)ds,  \label{dis_1}
\end{equation}
where $N^{\ast }$ is the decreasing rearrangement of $N$. Setting 
\begin{equation*}
s=\gamma
_{h}(B_{r})=N\omega _{N}\displaystyle\int_{0}^{r}e^{h(s)}s^{N-1}ds,
\end{equation*}
we get 
\begin{equation*}
\int_{0}^{\gamma _{h}(\Omega )}N^{\ast }(s)ds=N\omega
_{N}\int_{0}^{r^{\bigstar }}N^{\ast }(\gamma _{h}(B_{r}))e^{h(r)}r^{N-1}dr.
\end{equation*}
Note that 
\begin{equation*}
N^{\ast }(\gamma _{h}(B_{r}))=N(r),
\end{equation*}
since $N^{\ast }(\gamma _{h}(B_{r}))$ and $N(r)$ are equimeasurable and both
radially decreasing functions. Therefore 
\begin{equation}
N\omega _{N}\int_{0}^{r^{\bigstar }}N^{\ast }(\gamma
_{h}(B_{r}))e^{h(r)}r^{N-1}dr=N\omega _{N}\int_{0}^{r^{\bigstar
}}N(r)e^{h(r)}r^{N-1}dr=\int_{\Omega ^{\bigstar }}N(\left\vert x\right\vert
)d\gamma _{h}  \label{dis_2}
\end{equation}
Combining (\ref{dis_1}) and (\ref{dis_2}), we obtain the claim (\ref
{dis_BR_1}). Analogously it is possible to prove that 
\begin{equation}
\int_{\Omega }D(\left\vert x\right\vert )d\gamma _{h}\geq \int_{\Omega
^{\bigstar }}D(\left\vert x\right\vert )d\gamma _{h}.  \label{dis_BR_2}
\end{equation}
Indeed since $D$ is an increasing function, we have 
\begin{align*}
\int_{\Omega }D(\left\vert x\right\vert )e^{h(\left\vert x\right\vert )}dx&
\geq \int_{0}^{\gamma _{h}(\Omega )}D_{\ast }(s)ds \\
& =N\omega _{N}\int_{0}^{r^{\bigstar }}D_{\ast }(\gamma
_{h}(B_{r}))e^{h(r)}r^{N-1}dr=\int_{\Omega ^{\bigstar }}D(\left\vert
x\right\vert )d\gamma _{h},
\end{align*}
where $D_{\ast }$ is the increasing rearrangement of $D$. By (\ref{def G}),
(\ref{dis_BR_1}) and (\ref{dis_BR_2}), inequality (\ref{mu_1_i}) implies
\begin{equation*}
\mu _{1}(\Omega ;e^{h\left( |x|\right) })\leq \frac{\displaystyle 
\int_{\Omega ^{\bigstar }}\left( \left( w^{\prime }\left( \left\vert
x\right\vert \right) \right) ^{2}+\dfrac{N-1}{\left\vert x\right\vert ^{2}}
w\left( \left\vert x\right\vert \right) ^{2}\right) d\gamma _{h}}{ 
\displaystyle\int_{\Omega ^{\bigstar }}w(\left\vert x\right\vert
)^{2}d\gamma _{h}}=\mu _{1}(\Omega ^{\bigstar };e^{h\left( |x|\right) }),
\end{equation*}
which is our claim. Moreover, from the monotonicity properties
of the functions $N$ and $D$, it easy to realize that inequalities (\ref
{dis_BR_1}) and (\ref{dis_BR_2}) reduce to equalities only when $\Omega $ is
the ball $\Omega ^{\bigstar }.$

We finally exhibit an example showing that, in general, the condition about
the symmetry of the domain cannot be dropped.

\noindent Let 
\begin{equation*}
H_{n}(t):=(-1)^{n}e^{t^{2}}\left( \frac{d^{n}}{dt^{n}}e^{-t^{2}}\right) , 
\hspace{0.2cm}t\in \mathbb{R},
\end{equation*}
and 
\begin{equation*}
v_{n}(t):=H_{n}(t)e^{-t^{2}},\hspace{0.2cm}t\in \mathbb{R}.
\end{equation*}

Let $c$ and $d$ be the first and second positive zeros of $v_{5}^{\prime
}(t)=-8e^{-t^{2}}\left( 8t^{6}-60t^{4}+90t^{2}-15\right) $ respectively. It
is elementary to verify that 
\begin{equation}
c\in (0.43,0.44)\text{ and }d\in (1.33,1.34).  \label{c_d}
\end{equation}

We consider the following two-dimensional problem with anti-Gaussian
degeneracy 
\begin{equation*}
\left\{ 
\begin{array}{ll}
-\text{div}\left( e^{x^{2}+y^{2}}\nabla u\right) =\mu e^{x^{2}+y^{2}}u & 
\text{ in }T \\ 
&  \\ 
\dfrac{\partial u}{\partial \nu }=0 & \text{on }\partial T,
\end{array}
\right.
\end{equation*}
where $T=\left( c,d\right) \times \left( -c,c\right) $ . A straightforward
computation shows that $\mu _{1}(T)=12,$ moreover $\mu _{1}(T)$ is a double
eigenvalue and a corresponding set of independent eigenfunctions are $
u_{1}(x,y):=v_{5}(x)$ and $u_{2}(x,y):=v_{5}(y)$ (see, e.g., \cite{Taylor},
p.104 ff.).

Define 
\begin{equation*}
d\gamma _{2}=e^{x^{2}+y^{2}}dxdy,\text{ \ with }(x,y)\in \mathbb{R}^{2}.
\end{equation*}
Now we claim that the ball $B_{r_{T}}$, such that $\gamma
_{2}(B_{r_{T}})=\gamma _{2}(T)$, fulfills 
\begin{equation}
\mu _{1}\left( B_{r_{T}};\gamma _{2}\right) <12=\mu _{1}(T;\gamma _{2}).
\label{claim}
\end{equation}

\noindent Clearly 
\begin{equation*}
\gamma _{2}(B_{r})=\pi \left( e^{r^{2}}-1\right) =:\chi (r),\hspace{0.2cm} 
\text{with}\hspace{0.2cm}r>0,
\end{equation*}
and therefore 
\begin{equation*}
r_{T}=\chi ^{-1}\left( \gamma _{2}\left( T\right) \right) .
\end{equation*}
As recalled in Section 2, $\mu _{1}\left( B_{r_{T}},\gamma _{2}\right) $
satisfies the following variational characterization 
\begin{equation}
\mu _{1}\left( B_{r_{T}};\gamma _{2}\right) =\min \left\{ \frac{ 
\displaystyle \int_{B_{r_{T}}}\left\vert Dv\right\vert ^{2}d\gamma _{2}}{ 
\displaystyle \int_{B_{r_{T}}}
 v^{2}d\gamma _{2}} :
\text{ }v\in H^{1}\left( B_{r_{T}}\right) \backslash \left\{ 0\right\}
,\int_{B_{r_{T}}}vd\gamma _{2}=0\right\}  \label{B_T}
\end{equation}
In order to get an estimate from above for $\mu _{1}\left( B_{r_{T}};\gamma
_{2}\right) $, we use $v=x$ and $v=y$ as trial functions in (\ref{B_T})
obtaining 
\begin{equation*}
\mu _{1}\left( B_{r_{T}};\gamma _{2}\right) <\frac{\gamma _{2}(B_{r_{T}})}{ 
\displaystyle\int_{B_{r_{T}}}x^{2}d\gamma _{2}}\hspace{0.2cm}\text{ and } 
\hspace{0.2cm}\mu _{1}\left( B_{r_{T}};\gamma _{2}\right) <\frac{\gamma
_{2}(B_{r_{T}})}{\displaystyle\int_{B_{r_{T}}}y^{2}d\gamma _{2}}.
\end{equation*}
Summing up we get

\begin{equation}
\mu _{1}\left( B_{r_{T}};\gamma _{2}\right) <\frac{2\displaystyle 
\int_{0}^{r_{T}}e^{s^{2}}sds}{\displaystyle\int_{0}^{r_{T}}e^{s^{2}}s^{3}ds}
=k\left( r_{T}\right) ,  \label{mu1<k}
\end{equation}
where 
\begin{equation*}
k(r):=\frac{2e^{r^{2}}-2}{r^{2}e^{r^{2}}-e^{r^{2}}+1},\hspace{0.2cm}\text{
with}\hspace{0.2cm}r>0.
\end{equation*}
A Taylor expansion of $e^{x^{2}}$ allows to estimate from below $\gamma
_{2}(T)$ as follows

\begin{equation}  \label{gamma>2}
\gamma _{2}(T) > \int_{c}^{d} \left( 1+x^{2}+\frac{x^{4}}{2}+\frac{x^{6}}{6}
\right) dx \int_{-c}^{c} \left( 1+y^{2}+\frac{y^{4}}{2}+\frac{y^{6}}{6}
\right) dy>2,
\end{equation}
where the last inequality is an immediate consequence of (\ref{c_d}).

Since $\chi ^{-1}$ is an increasing function, by (\ref{gamma>2}) we have

\begin{equation*}
\chi ^{-1}\left( \gamma _{2}\left( T\right) \right) >\chi ^{-1}\left(
2\right) =\sqrt{\log \left( 1+\frac{2}{\pi }\right) }.
\end{equation*}

\noindent Finally since $k$ is a decreasing function, the above inequality
together with (\ref{mu1<k}) imply 
\begin{equation*}
\mu _{1}\left( B_{r_{T}};\gamma _{2}\right) < k\left( r_{T}\right) < k\left(
\chi ^{-1}\left(2\right) \right) =\frac{4}{(\pi +2)\log (1+\frac{2}{\pi })-2}
<12.
\end{equation*}
Hence the claim is proved.$\hfill \square $

\bigskip

\begin{remark}
Note that the assumption on the symmetry of $\Omega $ is used solely to
guarantee the orthogonality conditions \emph{(\ref{Orth_P_i})}.
\end{remark}

\bigskip

\noindent \textsc{Acknowledgements}. 
This paper was partially supported by the grants PRIN 2012 
``Elliptic and parabolic partial differential equations: 
geometric aspects, related inequalities, and applications'', 
FIRB 2013  ``Geometrical and qualitative aspects of PDE's'', 
and  STAR 2013 ``Sobolev-Poincar\'e inequalities: embedding constants, 
stability issues, nonlinear eigenvalues'' (SInECoSINE).


\bigskip

\begin{thebibliography}{99}
\bibitem{A} Ashbaugh M.S., \textsl{Isoperimetric and universal inequalities
for eigenvalues}, Spectral Theory and Geometry, Edinburgh, 1998, in: London
Math. Soc. Lecture Note Ser., vol. 273, Cambridge Univ. Press, Cambridge,
1999, 95--139.

\bibitem{AB} Ashbaugh M.S., Benguria R., \textsl{Sharp upper bound to the
first nonzero Neumann eigenvalue for bounded domains in spaces of constant
curvature}, J. Lond. Math. Soc. (2) 52 (2) (1995) 402--416.

\bibitem{Ba} Bandle C. \textsl{Isoperimetric inequalities and applications,}
Monographs and Studies in Mathematics 7, Pitman (Advanced Publishing
Program), Boston, Mass.-London, 1980.

\bibitem{BeLi} Benguria R.D., Linde H., \textsl{A second eigenvalue bound
for the Dirichlet Schr\"{o}dinger operator}, Comm. Math. Phys. 267 (3)
(2006) 741--755.

\bibitem{BBMP} Betta M. F., Brock F., Mercaldo A., Posteraro, M. R., \textsl{
\ \ \ Weighted isoperimetric inequalities on $\mathbb{R}^{n}$and
applications to rearrangements}, Math. Nachr. 281 (2008), no. 4, 466--498.

\bibitem{BCKT} Brandolini B., Chiacchio F., Krej\v{c}i\v{r}\'{\i}k D.,
Trombetti C., \textsl{The equality case in a Poincar\'{e}-Wirtinger type
inequality, } arXiv:1410.0676.

\bibitem{BCHT} Brandolini B., Chiacchio F., Henrot A., Trombetti C., \textsl{
\ \ \ Existence of minimizers for eigenvalues of the Dirichlet-Laplacian
with a drift,} arXiv:1406.6824.

\bibitem{BNT} Brasco L., Nitsch C., Trombetti C., \textsl{An inequality \`{a}
la Szeg\"{o}-Weinberger for the $p-$Laplacian on convex sets, }
arXiv:1407.7422.

\bibitem{BCM} Brock F., Chiacchio F., Mercaldo A., \textsl{Weighted
isoperimetric inequalities in cones and applications}, Nonlinear Anal. 75
(2012), no. 15, 5737--5755.

\bibitem{BMP} Brock F., Mercaldo A., Posteraro M. R., \textsl{On
isoperimetric inequalities with respect to infinite measures}, Rev. Mat.
Iberoam. 29 (2013), no. 2, 665--690.

\bibitem{cha1} Chavel I. \textsl{Lowest-eigenvalue inequalities}, Geometry
of the Laplace operator (Proc. Sympos. Pure Math., Univ. Hawaii, Honolulu,
Hawaii, 1979), pp. 79--89, Proc. Sympos. Pure Math., XXXVI, Amer. Math.
Soc., Providence, R.I., 1980.

\bibitem{Cha} Chavel I. \textsl{Eigenvalues in Riemannian Geometry}, New
York: Academic Press. 2001.

\bibitem{CdB} Chiacchio F., di Blasio G., \textsl{Isoperimetric inequalities
for the first Neumann eigenvalue in Gauss space, }Ann. I. H. Poincar\'{e} --
AN 29 (2012) 199--216.

\bibitem{CR} Chong K. M., Rice N. M., \textsl{Equimeasurable Rearrangements
of Functions}, Queen's Papers in Pure and Applied Mathematics, No. 28,
Queen's University, 1971.

\bibitem{C-H} Courant R., Hilbert D., \textsl{Methods of mathematical
physics }, vol. I and II, Interscience Publichers New York-London, 1966.

\bibitem{DPG1} Della Pietra F., Gavitone N., \textsl{Faber-Krahn Inequality
for Anisotropic Eigenvalue Problems with Robin Boundary Conditions}
Potential Analysis (2014) vol. 41 pag. 1147--1166.

\bibitem{DPG2} Della Pietra F., Gavitone N., \textsl{Stability results for
some fully nonlinear eigenvalue estimates } Communications in Contemporary
Mathematics Vol. 16, No. 5 (2014) 1350039 (23 pages).

\bibitem{H1} Henrot A., \textsl{Extremum problems for eigenvalues of
elliptic operators}, Frontiers in Mathematics. Birkh\"{a}user Verlag, Basel,
2006.

\bibitem{H2} Henrot A., Pierre M., \textsl{Variation et optimisation de
formes. Une analyse g\'{e}om\'{e}trique.}  Math\'{e}matiques  \&  Applications,
vol. 48, Springer, Berlin, 2005.

\bibitem{Ka} Kawohl B., \textsl{Rearrangements and Convexity of Level Sets
in PDE}, Lecture Notes in Mathematics 1150. New York: Springer Verlag, 1985.

\bibitem{Ke} Kesavan S., \textsl{Symmetrization \& applications}, Series in
Analysis, 3. World Scientific Publishing Co. Pte. Ltd., Hackensack, NJ, 2006.

\bibitem{KS} Kornhauser E.T. , Stakgold I., \textsl{A variational theorem
for\ } $\nabla ^{2}u+\lambda u =0 $ \textsl{and its applications,} J. Math.
Phys. 31 (1952) 45--54.

\bibitem{LS} Laugesen R.S., Siudeja B.A., \textsl{Maximizing Neumann
fundamental tones of triangles}, J. Math. Phys. 50 (11) (2009) 112903, 18 pp.

\bibitem{Mu} M\"{u}ller C., \textsl{Spherical Harmonics}, Lecture Notes
inMathematics, 17, Springer-Verlag, Berlin-New York 1966.

\bibitem{NS} Naito Y., Suzuki T., \textsl{Radial symmetry of self-similar
solutions for semilinear heat equations}, J. Differential Equations 163
(2000), no. 2, 407--428.

\bibitem{Rs} Rakotoson J. M., Simon B., \textsl{\ Relative rearrangement on
a measure space application to the regularuty of weighted monotone
rearrangement, I, II}, Appl. Math. Lett. 6 (1993), 75--78, 79--82.

\bibitem{S} Szeg\"{o} G., \textsl{Inequalities for certain eigenvalues of a
membrane of given area}, J. Rational Mech. Anal. 3 (1954) 343--356.

\bibitem{Taylor} Taylor M.E. , \textsl{Partial differential equations},
Vol.II, Qualitative Studies of Linear Equations. Appl. Math. Sciences 116,
Springer, N.Y. (1996).

\bibitem{W} Weinberger H.F., \textsl{An isoperimetric inequality for the
N-dimensional free membrane problem}, J. Rational Mech. Anal. 5 (1956)
633--636.
\end{thebibliography}
\end{document}